# Lattice-based Motion Planning for a General 2-trailer system

Oskar Ljungqvist[1], Niclas Evestedt[1], Marcello Cirillo[2], Daniel Axehill[1], Olov Holmer

*Abstract*— Motion planning for a general 2-trailer system poses a hard problem for any motion planning algorithm and previous methods have lacked any completeness or optimality guarantees. In this work we present a lattice-based motion planning framework for a general 2-trailer system that is resolution complete and resolution optimal. The solution will satisfy both differential and obstacle imposed constraints and is intended as a driver support system to automatically plan complicated maneuvers in backward and forward motion. The proposed framework relies on a precomputing step that is performed offline to generate a finite set of kinematically feasible motion primitives. These motion primitives are then used to create a regular state lattice that can be searched for a solution using standard graph-search algorithms. To make this graph-search problem tractable for real-time applications a novel parametrization of the reachable state space is proposed where each motion primitive moves the system from and to a selected set of circular equilibrium configurations. The approach is evaluated over three different scenarios and impressive real-time performance is achieved.

## I. INTRODUCTION

In this paper we present a resolution complete and resolution optimal motion planning framework for a general 2-trailer system in both forward and backward motion. The general 2-trailer system is highly nonlinear and even unstable in backward motion which makes path planning for this system challenging. To enable the effective use of graph search algorithms for path planning under the kinematic constraints imposed by this system, a novel parametrization of the state lattice is proposed. By calculating motion primitives that move the system from and to a chosen set of equilibrium configurations, two system states can be directly removed from the state lattice and make real-time use of classical graph search algorithms tractable. The motion planner could be used as a driver support system to plan complex maneuvers in parking scenarios or as a stand alone planner for autonomous maneuvering with trailers. To the best of the author's knowledge this work presents the first resolution complete motion planning framework for a reversing general 2-trailer system. The focus of this paper is to generate paths that are kinematically feasible. Therefore, the reader is referred to our previous work on path following controllers in order to find techniques to stabilize the system around these paths [1], [2].

### A. Related work

The nonlinear dynamics of a standard trailer configuration with the hitch connection in the center of the rear axle are well understood and the derivation of the equations for a standard n-trailer configuration can be found in [3]. For the standard 1-trailer configuration an exact local motion planner was presented in [4]. In [5] it is shown that $\eta^4$-splines can be used to generate smooth and feasible paths for the standard 1-trailer system and these splines are used in a local planner. Trajectory generation for the standard n-trailer is studied in [6] where the kinematic equations are first transformed into chained form and steered from an initial configuration to a goal configuration before a transformation back to the original coordinates is performed. One and two trailer planning simulations are presented but obstacles are omitted. Impressive results on global planning with obstacles for a standard 2-trailer is presented in [7]. However, the approaches presented above consider the case with on-axle hitching despite that most practical applications have off-axle hitching making the kinematics more challenging and the system equations more complicated.

To include the general 2-trailer system within a motion planning framework, we presented a probabalistic sampling based approach in [8]. Even though the planner was capable of solving several complicated problems, the framework lacks any guarantees for completeness and optimality. To overcome the lack of guarantees [9] introduced lattice-based planners which offer resolution completeness and optimality. By discretizing the state space of the model and pre-compute a set of feasible motions that moves the system from one discretization point to another, a state lattice can be formed. In this way the differential constraints have already been considered offline and during planning the planner only needs to perform a search over the available motions.

Lattice-based planners have been deployed with great success on several robotic platforms [9], [10], [11]. However, a problem with lattice-based approaches is the exponential complexity in the dimension of the state space which can limit the use for more complicated models.

In this work we have circumvented this problem by parametrizing the state lattice such that the motions always moves to system from and to a circular equilibrium configuration. In this way two states can be directly removed from consideration during planning which makes the use of lattice-based planners tractable for the general 2-trailer system. By only planning between equilibrium configurations the movements are more restricted but we show through several example scenarios that the degrees of freedom for the system is enough to complete several challenging and relevant scenarios. The generation of the motion primitives is performed using a solver for the Two-Point Boundary Value Problem (TPBVP) of the general 2-trailer system that was developed in [12]. The solver is based on the ACADO optimization toolkit [13] and is used to create both forward and backward motions. To overcome the stability problem when generating motion primitives for backward motion, a symmetry result is established which shows that it is kinematically feasible to follow a path in backward motion that was generated by driving the model forward. This symmetry property can then be exploited to efficiently generate the motion primitives for backward motion in a numerically stable way.

*The research leading to these results has been carried out within the iQMatic project funded by FFI/VINNOVA.
[1]Division of Automatic Control, Linköping University, Sweden, (e-mail: {niclas.evestedt, oskar.ljungqvist, daniel.axehill}@liu.se)
[2]Scania Technical Centre, Södertälje Sweden, (e-mail: marcello.cirillo@scania.com)

The outline of the remainder of the paper is as follows: In Section II a kinematic model of the general 2-trailer system and a symmetry result for a certain class of driftless systems is presented. In Section III the lattice-based motion planner is presented and Section IV describes the framework that is used to generate the set of kinematically feasible motion primitives. Finally, Section V and Section VI present the results and conclusions, respectively.

## II. System Model

The kinematic model of a general 2-trailer system with off-axle hitching on the pulling vehicle is presented in [14]. The vehicle configuration is schematically illustrated in Fig. 1. To model this nonholonomic system the generalized coordinates $\mathbf{p} = (x_3, y_3, \theta_3, \beta_3, \beta_2)$ are used, where $(x_3, y_3)$ denotes the center of the rear axle of the trailer, $\theta_3$ is the global orientation of the trailer, $\beta_3$ is the relative angle between the trailer and the dolly and $\beta_2$ is the relative angle between the dolly and the truck. The geometric lengths $L_3$, $L_2$ and $L_1$ are the distance between the rear axle of the trailer to the rear axle of the dolly, the distance between the rear axle of the dolly to the off-axle hitch connection of the truck and the distance between the axles of the truck, respectively. The off-axle hitch length of the truck is denoted $M_1$. The kinematic model for this general 2-trailer becomes

$$\dot{x}_3 = v\cos\beta_3 \cos\beta_2 \left(1 + \frac{M_1}{L_1}\tan\beta_2 \tan\alpha\right)\cos\theta_3 \quad \text{(1a)}$$

$$\dot{y}_3 = v\cos\beta_3 \cos\beta_2 \left(1 + \frac{M_1}{L_1}\tan\beta_2 \tan\alpha\right)\sin\theta_3 \quad \text{(1b)}$$

$$\dot{\theta}_3 = v\frac{\sin\beta_3 \cos\beta_2}{L_3}\left(1 + \frac{M_1}{L_1}\tan\beta_2 \tan\alpha\right) \quad \text{(1c)}$$

$$\dot{\beta}_3 = v\cos\beta_2 \left(\frac{1}{L_2}\left(\tan\beta_2 - \frac{M_1}{L_1}\tan\alpha\right) - \frac{\sin\beta_3}{L_3}\left(1 + \frac{M_1}{L_1}\tan\beta_2 \tan\alpha\right)\right) \quad \text{(1d)}$$

$$\dot{\beta}_2 = v\left(\frac{\tan\alpha}{L_1} - \frac{\sin\beta_2}{L_2} + \frac{M_1}{L_1 L_2}\cos\beta_2 \tan\alpha\right) \quad \text{(1e)}$$

where $v$ is the longitudinal velocity for the rear axle of the truck and $\alpha$ is the steering angle of the truck. The model is valid under a no-slip assumption and is also assumed to operate on a relatively flat surface. Since the operational speed for maneuvering with a general 2-trailer is quite low these assumptions are expected to hold. The direction of motion is important for the stability of the system (1), in forward motion ($v > 0$) the system is marginally stable and in backward motion ($v < 0$) the system is unstable and the truck and trailer angles can fold and enter what is called a jack-knife state.

Since the velocity enters linearly in (1), a method known as time-scaling can be applied to eliminate the longitudinal speed dependence from the model. This means that the velocity profile can be constructed separately from the path construction. Therefore, we only consider the longitudinal speed $v$ to take on the values $v = 1$ for forward motion and $v = -1$ for backward motion when generating the motion primitives. In practice, we have saturation in the steering angle $|\alpha| \leq \alpha_{\max}$ and on the rate of the steering angle $|\omega| \leq \omega_{\max}$. These constraints have to be considered during the motion primitive generation such that the constructed motions are feasible to follow for a real system.

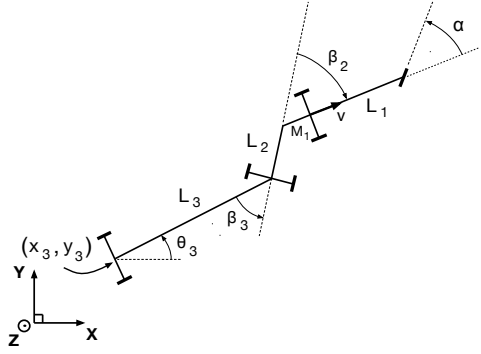

Fig. 1: Schematic view of the configuration used to model the general two-trailer system used in this work.

### A. Symmetry

A system is symmetric if a state trajectory that brings the system from an initial state to a final state can be reversed in time and the kinematic constraints will also hold for the time reversed trajectory going from the final state to the initial state. For a class of driftless systems, including the model for the general 2-trailer presented in (1), that can be written on the form

$$\dot{\mathbf{p}}(t) = v(t)f(\mathbf{p}(t), u(t)) \quad \text{(2)}$$

where $\mathbf{p}$ denotes the states, $v(t)$ and $u(t)$ denotes the control signals, the symmetry property can be shown to hold. The result is summarized in Theorem 1.

*Theorem 1:* The driftless system (2) is symmetric when the time-reversing control signals:

$$\bar{v}(t) = -v(T-t), \quad \bar{u}(t) = u(T-t)$$

are applied and the result is a reversed state trajectory $\bar{\mathbf{p}} = \mathbf{p}(T-t)$, $t \in [0, T]$.

*Proof:* If we assume $\bar{\mathbf{p}}(t) = \mathbf{p}(T-t)$, $t \in [0, T]$ then

$$\frac{d}{dt}\bar{\mathbf{p}}(t) = \frac{d}{dt}\mathbf{p}(T-t) = \{\tau = T-t\}$$
$$= \frac{d\tau}{dt}\frac{d}{d\tau}\mathbf{p}(\tau)$$
$$= -v(\tau)f(\mathbf{p}(\tau), u(\tau))$$
$$= \bar{v}(t)f(\bar{\mathbf{p}}(t), \bar{u}(t))$$

This together with the initial condition $\bar{\mathbf{p}}(0) = \mathbf{p}(T)$ concludes the proof. ∎

The symmetry property will be used in the motion primitive generation in Section IV to avoid that the unstable dynamics in backward motion can cause numerical problems while solving the TPBVP. Instead the problem is solved from the final state to the initial state in forward motion and the resulting control signal and state trajectory can then be time-reversed by using the result from Theorem 1 and applied to generate the backward motion segment.

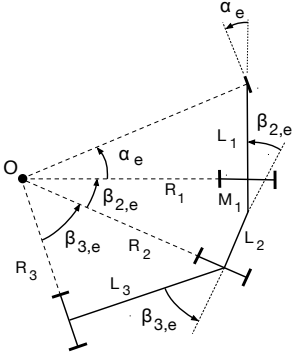

Fig. 2: Stationary equilibrium point for steering angle $\alpha_e$. The system will travel in a circular path with a radius determined by the geometry and $\alpha_e$.

### B. Equilibrium configurations

In order to reduce the dimension of the state lattice the motion primitives will always be calculated between two equilibrium configurations of the system. Given a constant steering angle $\alpha_e$ there exists a circular equilibrium configuration, as depicted in Fig. 2, where $\dot{\beta}_3$ and $\dot{\beta}_2$ are equal to zero. At this equilibrium the truck and trailer system will travel along circular arcs with radiuses determined by $\alpha_e$. The equilibrium configuration can be determined using trigonometry [2] which gives the following relations:

$$\beta_{3e} = \text{sign}(\alpha_e)\arctan\left(\frac{L_3}{R_3}\right) \qquad (3a)$$

$$\beta_{2e} = \text{sign}(\alpha_e)\left(\arctan\left(\frac{M_1}{R_1}\right) + \arctan\left(\frac{L_2}{R_2}\right)\right) \qquad (3b)$$

where $R_1 = L_1/|\tan\alpha_e|$, $R_2 = \sqrt{R_1^2 + M_1^2 - L_2^2}$, $R_3 = \sqrt{R_2^2 - L_3^2}$. Now, instead of discretizing the full state-space of the model during the lattice creation, $\beta_3$ and $\beta_2$ will be determined from the given equilibrium configuration in the initial and goal state. This directly reduces two states from the dimension of the search space and only $\alpha_e$ needs to be considered in the lattice creation.

### III. LATTICE PLANNER

Given a model of vehicle kinematics, the intuition behind lattice-based motion planning is to sample the state space of the model in a regular fashion and constrain the motions of the vehicle to a lattice graph $\mathcal{G} = \langle \mathcal{V}, \mathcal{E} \rangle$, that is, a graph embedded in a Euclidean space $\mathbb{R}^n$ which forms a regular tiling [9], [15]. Each vertex $v \in \mathcal{V}$ represents a discrete state of the vehicle while each edge $e \in \mathcal{E}$ encodes a motion which respects the nonholonomic and physical constraints of the vehicle.

The discretization of the lattice defines which states the vehicle can reach and the constraints of the vehicle is encoded in the motion primitive set $P$. At a valid state the kinematic model of the general 2-trailer (1) can be represented by a four-dimensional vector $\mathbf{s} = (x_3^d, y_3^d, \theta_3^d, \alpha_e^d)$, where $d$ highlights a discretized state. The positions $(x_3^d, y_3^d)$ lies on a grid of resolution $r$, $\theta_3^d \in \Theta$ and $\alpha_e^d \in \Phi$, where $\Theta$ and $\Phi$ are a finite set of allowed orientations of the trailer and of allowed equilibrium steering angles, respectively. Under some assumptions the model of the general 2-trailer system (1) is position-invariant and thus we can design $P$ to be position-invariant. Every motion primitive $p \in P$ is calculated by using a TPBVP solver to connect a set of initial states $\mathbf{s}_i = (0, 0, \theta_{3,i}^d, \alpha_{e,i}^d)$ to a set of neighbouring states in a discrete, bounded neighborhood in free space. The TPBVP solver guarantees that the motions respect the kinematic and physical constraints of the vehicle, while the position-invariant property ensures that the primitives are translatable to other states. The method used to generate the set of motion primitives $P$ will be further explained in Section IV.

Finally, a cost $g(p)$ is associated with each $p \in P$. In our implementation, $g(p)$ is calculated by multiplying the distance covered by the rear axle of the trailer by a cost factor which penalizes backward motion and an additional constant cost that is added for every direction change.

A planning problem is defined by an initial state $\mathbf{s}_I$, a goal state $\mathbf{s}_G$ and a world representation $\mathcal{W}$, in which all known obstacles are included. A *feasible solution* is a sequence of collision-free primitives $(p_0, \ldots, p_n)$ connecting $\mathbf{s}_I$ to $\mathbf{s}_G$. Given the set of all feasible solutions to a problem, the optimal solution is the one with minimum cost. Here, we explore the state space using the standard $A^*$ search algorithm where an admissable free-space heuristic table, as described in [15], is precomputed and used for exploration. By using this search method both resolution completeness and resolution optimality can be guaranteed.

### IV. LATTICE CREATION

In order to generate a state lattice, the state space of the model (1) needs to be discretized. In this work the position has a discretization accuracy $r = 0.5$ m for both $x_3^d$ and $y_3^d$ and the orientation $\theta_3^d$ is discretized into $|\Theta| = 16$ different different angles. The equilibrium steering angle $\alpha_e^d$ is discretized into $|\Phi| = 3$ different angles, $\alpha_e^d \in \Phi = \{-0.2117, 0, 0.2117\}$, where $|\alpha_e^d| = 0.2117$ corresponds to circular movements according to (3) with $R_3 = 20$ m for our full-size test vehicles. To ensure smooth transitions between every motion segment within the lattice the motion between the initial and the final state must be solved between the equilibrium configuration specified by the steering angle in every grid point. In order to guarantee that the generated motions are feasible to execute for a real system first and second order derivatives of the steering angle $\alpha$ are added to the model. This is done to be able to ensure that the steering angle becomes a $C^1$-function in time and physical constraints on the steering angle rate $\omega$ and acceleration $u$ can be considered. The steering angle acceleration is then used as control signal $u$ in the optimal control problem (4). Introduce the new state vector $\bar{\mathbf{p}} = (\mathbf{p}, \alpha, \omega)$. To generate

motion segments between the grid points the following TPBVP is solved:

$$\begin{aligned}
\underset{u(\cdot),\,T}{\text{minimize}} \quad & \int_0^T f_0(\bar{\mathbf{p}}(t), u(t))\, dt \quad (4)\\
\text{subject to} \quad & \dot{\mathbf{p}}(t) = f(\mathbf{p}(t), \alpha(t)),\\
& \dot{\alpha}(t) = \omega,\ \dot{\omega}(t) = u,\\
& \bar{\mathbf{p}}(0) = (\mathbf{p}_i^d, \alpha_{e,i}^d, 0),\ \bar{\mathbf{p}}(T) = (\mathbf{p}_f^d, \alpha_{e,f}^d, 0),\\
& |\beta_3(t)| \le \beta_{3,\max},\ |\beta_2(t)| \le \beta_{2,\max},\\
& |\alpha(t)| \le 0.8\alpha_{\max},\ |\omega(t)| \le \omega_{\max},\\
& |u(t)| \le u_{\max}
\end{aligned}$$

where $f_0(\bar{\mathbf{p}}(t), u(t))$ is the objective function and $\dot{\mathbf{p}}(t) = f(\mathbf{p}(t), \alpha(t))$ is the model of the general 2-trailer system described in (1). The objective function used to generate the motion primitives was chosen to $f_0(\bar{\mathbf{p}}(t), u(t)) = 10\omega^2 + u^2$ which make it awarding to generate smooth steering angles when possible. The constraints $\bar{\mathbf{p}}(0) = (\mathbf{p}_i^d, \alpha_{e,i}^d, 0)$ and $\bar{\mathbf{p}}(T) = (\mathbf{p}_f^d, \alpha_{e,f}^d, 0)$ are the initial and final state conditions determined from the grid, where the steering angle rate is constrained to be zero at the initial and the final state to ensure that the steering angle is a $C^1$-function in time, even when different motion segments are combined in the lattice planner. Physical constraints are added for the trailer angles, the steering angle, the steering angle rate and acceleration with maximum absolute values $\beta_{3,\max} = \pi/2$, $\beta_{2,\max} = \pi/2$, $\alpha_{\max} = \pi/4$, $\omega_{\max} = 1.5$ and $u_{\max} = 40$, respectively. Additionally, in order to leave room for a low-level path-following controller [1] to be able to reject disturbances during the path execution, the constraints on the steering angle are additionally 20% tightened such that the generated primitives do not saturate the maximum physical steering angle. The optimal control problem in (4) can now be solved using numerical optimal control. In this work the ACADO Toolkit is used [13].

### A. Generation of Motion Primitives

Even though the motion primitive generation is done offline it is not desirable to make an exhaustive generation of primitives to all grid points due to computation time and the high risk of creating redundant segments. Instead, a careful selection of configurations has been performed and primitives are only computed to grid points within this selected set. The same set of points are used for both forward and backward motions but the symmetry result established in Theorem 1 is used for generation of backward motion segments. In this work every transition between two equilibrium configurations is calculated with three different complexity levels which each requires different amount of steering input to manage the transition. Fig. 3 illustrates these complexity levels for a reversing transition from an initial configuration with $(\theta_{3,i}^d, \alpha_{e,i}^d) = (0, 0)$ to a final configuration with $(\theta_{3,f}^d, \alpha_{e,f}^d) = (\pi/2, 0)$. By changing the final position $(x_3^d, y_3^d)$ different amounts of maneuvering space is given. The closer the final position is to the initial position the more steering angle effort is required. This idea is then repeated for all possible initial configurations to different final configurations and as an example the complete set of possible motions from $\theta_{3,i}^d = 0$ is shown in Fig. 4. By rotating this set into all possible discretization points for the heading the full set of motions can be created.

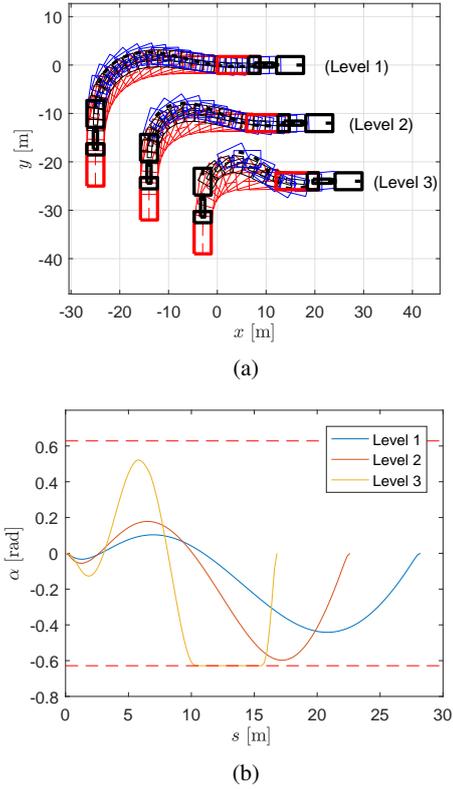

(a)

(b)

Fig. 3: An illustration of the three complexity levels that are used to generate the motion primitives when performing a transition from $(\theta_{3,i}^d, \alpha_{e,i}^d) = (0, 0)$ to $(\theta_{3,f}^d, \alpha_{e,f}^d) = (\pi/2, 0)$ in backward motion. In (a) the traversed paths are plotted and in (b) the steering angles are plotted against the travelled distance of the path taken by the rear axle of the trailer for the three levels where the red dotted lines correspond to the bounds on the steering angle. As can be seen the steering angle effort increases significantly when narrowing the maneuvering distance.

### B. Motion Primitive Reduction

To further improve the search speed the inital primitive set $P$ is reduced using the techniques described in [16]: given an initial configuration, all redundant primitives $p$ are removed from $P$ if the state transition for $p$ can be achieved using a combination of other primitives already in the set and the new combined cost $c_{\text{comb}} \le c_p m$, where $m \ge 1$.

Intuitively, we remove all primitives from the set $P$ that can be substituted by a combination of other primitives if the cost of the alternative combination is at most $m$ times worse than the original. In our tests, the original primitives were reduced using $m = 1.2$ which resulted in a new set $P'$ with a reduction factor of about 20 %. For each initial state in $P$ the number of applicable primitives ranged from 64 to 128 with a total number $|P| = 4096$ primitives while in the reduced set $P'$ the range was between 57 and 92 and a total number of $|P'| = 3296$ primitives.

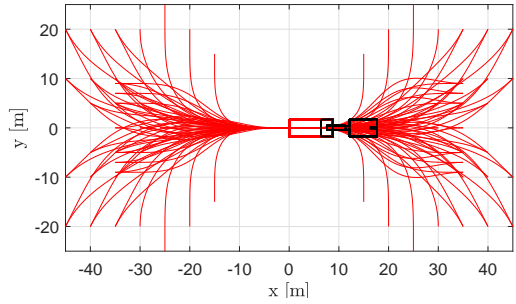

(a) The set of motion primitives for $(\theta_{3,i}^d, \alpha_{e,i}^d) = (0,0)$ to different final states on the lattice grid.

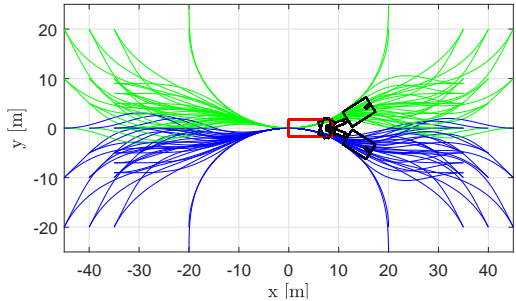

(b) The set of motion primitives for $(\theta_{3,i}^d, \alpha_{e,i}^d) = (0, 0.2117)$ (green) and $(\theta_{3,i}^d, \alpha_{e,i}^d) = (0, -0.2117)$ (blue) to different final states on the lattice grid.

Fig. 4: The set of motion primitives for zero initial heading of the trailer and different initial equilibrium configuration states to different final states on the lattice grid. The colored paths are the paths taken by the rear axle of the trailer $(x_3, y_3)$ during the different motions.

## V. RESULTS

The performance of the lattice planner is evaluated through a series of simulations performed on a standard laptop with an Intel Core i7-4810QM CPU @ 2.80GHz. In these examples, we run standard $A^*$ searches to evaluate the time required to find a solution. For the reduced primitive set $P'$, we calculated a free-space heuristic table, as described in [15]. The motion primitives are generated with lengths corresponding to our full-size test vehicles with $L_3 = 7.59$ m, $L_2 = 3.75$ m, $L_1 = 4.66$ m and $M_1 = 0.8$ m. In all scenarios the goal state $\mathbf{s}_G$ of the trailer is marked with a blue box where the white arrow represents the orientation. At the goal state the equilibrium steering angle is constrained to zero. The world $\mathcal{W}$ is represented as a grid-map and the regions which are occupied with obstacles are marked in red.

### A. Open Area

As a first scenario an open area planning problem is constructed. The objective is to find a solution that moves the truck and trailer system about 5 meters in the lateral direction and change the orientation of the trailer with an angle of 180 degrees. The scenario is shown in Fig. 5 together with the solution calculated by the planner. The planning time was only 76 milliseconds and the search tree expanded 233 edges. The white path represent the planned path for the rear axle of the trailer and the colored points represent the path to the first direction change. As can be seen the solution is mainly constructed by combining three types of motion primitives; A parallel movement and two 90 degrees turns in backward and forward motion, respectively.

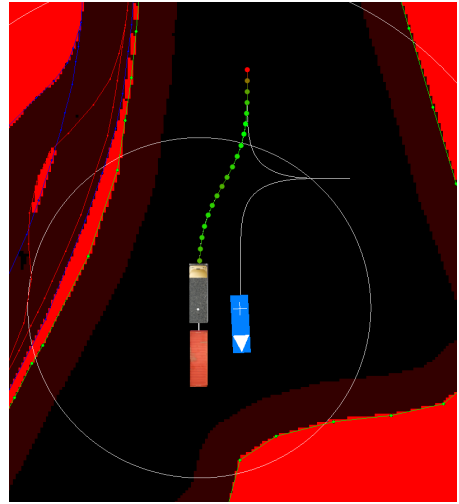

Fig. 5: Open area scenario where the vehicle configuration has to be turned around by combining forward and backward motion. The white and colored path is the planned path taken by the rear axle of the trailer during this maneuver.

### B. Reverse Parking

A common scenario a truck driver encounters is to navigate a truck and trailer in a parking lot and park the trailer by reversing into a free parking slot. The intention with the maneuver is either to park and disconnect the trailer or unload its cargo. The scenario setup is shown in Fig. 6 where the vehicle must navigate to an appropriate position from where it is possible to perform a final reversing maneuver into the free parking slot. In this scenario the planning time was 39 milliseconds and the search tree expanded only 102 edges. In this scenario the solution is constructed by combining motion primitives that performs 90 degrees turns, parallel movements with a small change in orientation and short forward and backward segments.

### C. Parking Lot

The last scenario is a planning problem where the objective is to find a path out from a parking lot. The parking lot exposes the planner to a complicated environment. The scenario can be seen in Fig. 7. By combining forward and backward motion segments the planner finds the solution in 161 milliseconds and the $A^*$ search algorithm expanded 454 edges. A drawback of lattice based approaches can be seen in the figure where the trailer orientation does not end up in the selected final orientation at the goal. This is due to the discretization in the lattice and if the desired goal state does not perfectly fit on a grid point the closest state in the lattice must be chosen as the goal state. This can be alleviated by the use of numerical optimal control as a post-processing step to end up in the goal state exactly.

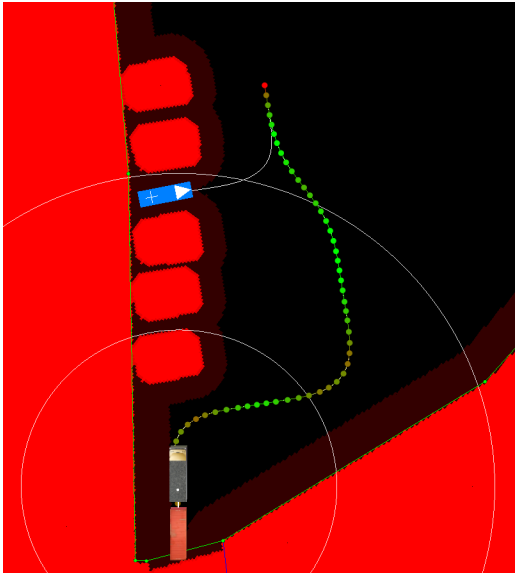

Fig. 6: Reverse parking scenario where the goal is to park the trailer in the free parking lot (the blue box). The white and colored path is the planned path for the rear axle of the trailer during this maneuver.

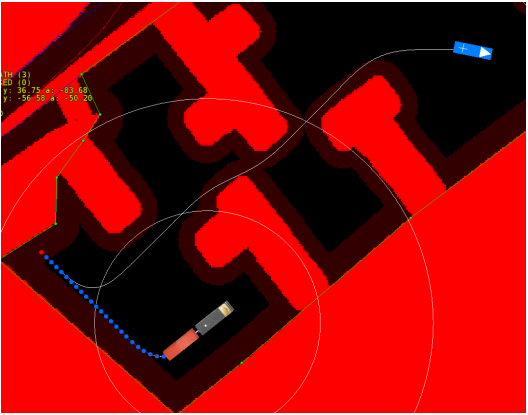

Fig. 7: Parking lot scenario where the vehicle should find its way out of the parking lot. The white and colored path is the path taken by the rear axle of the trailer during this maneuver.

## VI. CONCLUSIONS AND FUTURE WORK

In this paper a lattice-based motion planning framework for a general 2-trailer vehicle configuration is presented. A novel method for generating motion primitives is established which enables the user to design the motions such that they are both kinematically feasible and also satisfy physical constraints on the vehicles. A symmetry result is established for a certain class of driftless systems, including the general 2-trailer system, which is used to eficently generate motion primitives in backward motion. The generated motion primitives are then used within a lattice-based planning framework which offer guarantees on both resolution completeness and resolution optimality. To enable real-time performance of the lattice planner the state space is parametrized such that the motions always traverse the system from and to a circular equilibrium configuration. The motion planner is evaluated over three different scenarios where the planner is able to find a solution within fractions of a second. However a big drawback is the discretization of the state space which limits manuverability and prevents goal positions outside the discretized state space to be reached exacly. This problem could be handled by a post-processing step using numerical optimal control which brings the solution to the goal position exactly. Furture and ongoing work, include the integration of the planner on a full-sized truck and trailer together with the path following controller presented in [1].